\documentclass[10pt]{article}
\usepackage{amsmath,amssymb,amsfonts,oldlfont,graphics,latexsym,dsfont}
\usepackage{graphics}
\usepackage{graphicx}
\usepackage{epic,tikz}
\usepackage{geometry,verbatim}
\usetikzlibrary{calc}
\usetikzlibrary{arrows}

\usepackage{algorithm}
\usepackage{algorithmic}

\newcommand{\cp}{\ensuremath{\,\Box\,}}

\newcommand{\proof}{\noindent{\bf Proof\, }}
\newcommand{\qed}{\hfill $\square$\bigskip}

\newcommand{\bigbox}{\mathop{\raisebox{-0.25ex}{\hbox{\Large{$\square\hspace{0.5pt}\vspace{0.5pt}$}}}}}

\newtheorem{theorem}{Theorem}
\newtheorem{lemma}[theorem]{Lemma}

\newtheorem{definition}{Definition}

\newcommand{\A}{\ensuremath{\operatorname{Aut}}}

\title{Cartesian products of directed  graphs with loops}

\author{
Wilfried Imrich\thanks{Partially supported by OEAD Projekt SI 08/2016.}\\
\small Montanuniversit\"at Leoben\\
\small Franz Josef-Stra\ss e 18, 8700 Leoben, Austria\\
\small\tt imrich@unileoben.ac.at\\
\and Iztok Peterin\thanks{Partially supported by the ARRS under the research
grant P1-0297.}\\
\small University of Maribor, FEECS\\
\small Smetanova 17, 2000 Maribor, Slovenia\\
\small\tt iztok.peterin@um.si
\date {}}

\begin{document}

\maketitle

\begin{abstract} We show that every nontrivial finite or infinite connected directed graph  with loops and at least one  vertex without a loop is uniquely representable as a Cartesian or weak Cartesian product of prime graphs. For finite graphs the factorization can be computed in linear time and space.
\end{abstract}


\bigskip


\section{Introduction}

This note treats finite and infinite directed graphs with or without loops. It is shown that every connected, finite or infinite directed graph with at least one vertex without loop is uniquely representable as a Cartesian or weak Cartesian product of prime graphs, and that the factorization can be computed in linear time and space for finite  graphs.

The note extends and unifies results by Boiko et al.~\cite{BoCu-2016} about the Cartesian product of finite undirected graphs with loops, and   by Crespelle and Thierry~\cite{CreTh-2015} about finite directed graphs. For infinite graphs it generalises  a result by Miller \cite{miller-1970} and Imrich \cite{im-1971} about the weak Cartesian product.

Let us briefly mention that unique prime factorization with respect to the Cartesian product of connected
finite graphs was first shown 1960 by Sabidussi\footnote{1963 an   independent proof was published by Vizing \cite{vizing-1963}.} \cite{sabi-1960}, and that Sabidussi also introduced the weak Cartesian product.

Sabidussi's  proof is non-algorithmic. For undirected graphs the first  factorization algorithm  is due to Feigenbaum, Hershberger and Sch\"{a}ffer \cite{fehe-1985}. Its complexity  is $O(n^{4.5})$, where $n$ is the number of vertices of the graph. Subsequently the  complexity was further reduced by a number of authors.  The latest improvement, Imrich and Peterin  \cite{impe-2007}, is  linear in time and space in the number of edges.

For directed graphs the first factorization  algorithm  is from Feigenbaum \cite{fe-1986} and assumes the undirected decomposition provided by \cite{fehe-1985}. Crespelle and Thierry \cite{CreTh-2015} also assume an undirected decomposition and then compute the prime factorization of the directed graph in linear time and space.  Here we present a considerably simpler algorithm of the same complexity and extend it to the case when loops are allowed. We use the same data structure as in  \cite{impe-2007} and wish to  remark  that  a slight variation of the algorithm in \cite{impe-2007} would also allow a direct computation of the prime factors of connected directed graphs (with or without loops)  in linear time.\bigskip

\section{Preliminaries}
\label{sec:prelim}

A directed graph $G$ with loops consists of a set $V(G)$ of vertices together with a subset $A(G)$ of $V(G)\times V(G).$ The elements of $A(G)$ are called \emph{arcs} and are ordered pairs of vertices. If $ab$ is an arc, we call $a$ its \emph{origin}, $b$ its \emph{terminus}, and also refer to $a$ and $b$ as \emph{endpoints}.

We allow that   $a$ equals $b.$ In this case we speak of a  \emph{loop} at vertex $a$ and say the vertex $a$ is \emph{looped}.  Notice that it is possible that $A(G)$ contains both $ab$ and $ba.$

To make the notation better readable we often write $v\in G$ instead of $v\in V(G)$, and $e\in G$  instead of
$e \in A(G)$ or $E(G).$

We introduce the notation $\overrightarrow{\Gamma_0}$ for the class of directed graphs with loops.
To every graph $G \in \overrightarrow{\Gamma_0}$ we also define its \emph{shadow} ${\mathcal{S}}(G).$ It has the same vertex set as $G$ and its set of edges $E(G)$ consists of all unordered pairs $\{a, b\}$ of distinct vertices for which $ab$, $ba$ or both are in $A(G).$ To indicate that $\{a, b\}$ is an edge, we will use the notation $[a,b]$, or simply $ab$.
We say $G$ is \emph{connected} if  ${\mathcal{S}}(G)$ is connected, and set the\emph{ distance} $d_G(u,v)$ between two vertices, and the\emph{ degree} $d_G(u)$ of a vertex,  equal to  $d_{{\mathcal{S}}(G)}(u,v)$, resp. $d_{{\mathcal{S}}(G)}(u)$. The minimum degree of $G$ is denoted by $\delta$.
We also use the notation   $\Gamma$ for the class of simple graphs,   $\Gamma_0$ for the class of simple graphs with loops, and $\overrightarrow{\Gamma}$ for  the class of directed graphs without loops.

The \emph{Cartesian product $G\cp H$ of graphs in} $\overrightarrow{\Gamma_0}$ is defined on the Cartesian product $V(G)\times V(H)$ of the vertex sets of the factors. Its set of arcs  is
$$A(G\cp H)  =  \big\{ (x,u)(y,v) \mid xy\in A(G) \text{ and } u=v,  \text{ or, }  x=y \text{ and } uv\in A(H))\big\}.$$

If $G$ and $H$ have no loops, then this is also the case for $G\cp H.$ To obtain the definition of the Cartesian product of undirected graphs, one just replaces $A(G)$ by $E(G).$ Hence, the new definition generalizes the definition of the Cartesian product of simple graphs, directed graphs, and simple graphs with loops. Note that  ${\mathcal{S}}(G\cp H) = {\mathcal{S}}(G)\cp {\mathcal{S}}(H).$

Clearly Cartesian multiplication is commutative and the trivial graph $K_1$ is a unit.  It is well known that it is associative in $\Gamma.$ That it is associative in  $\overrightarrow{\Gamma}$ and $\Gamma_0$ was shown in \cite{fe-1986}, resp.  \cite{BoCu-2016}, for finite graphs.  We defer the proof that this also holds for finite graphs in $\overrightarrow{\Gamma}_0$ to   Section \ref{sec:weak}, where  we show  associativity in $\overrightarrow{\Gamma}_0$ with respect to the Cartesian and the weak Cartesian product. 

A nontrivial, connected graph $G$ with at least one unlooped vertex is called  \emph{irreducible} or \emph{prime} with respect to  Cartesian multiplication if, for every factorization $G= A \cp B$, either $A$ or $B$ has only one vertex.
Every finite  connected  graph in $\overrightarrow{\Gamma}_0$ with at least one unlooped vertex is uniquely representable as a  Cartesian product, up to the order and isomorphisms of the factors. Again we defer the proof and show in Section \ref{sec:unique} that unique factorization holds for finite graphs with respect to the Cartesian product and for infinite graphs with respect to the Cartesian or the weak Cartesian product.

For our proofs and algorithms projections and layers play an important role. The  \emph{$i^{\text{th}}$ projection} $p_i:V(G)\rightarrow V(G_i)$  of a product $\prod_{i=1}^kG_i$ is defined by $(v_1,\ldots,v_k)\mapsto v_i$ and the \emph{$G_i$-layer $G_i^v$  through a vertex} $v\in G$ is the subgraph induced by the set
$$\{w\in V(G)\,|\, w_j = v_j \text{ for all } j\neq i, 1\leq j \leq k\}.$$
Sometimes we will also use the notation $p_{G_i}$ instead of $p_i.$

\begin{figure}[htb]
\begin{center}
\includegraphics[width=11cm]{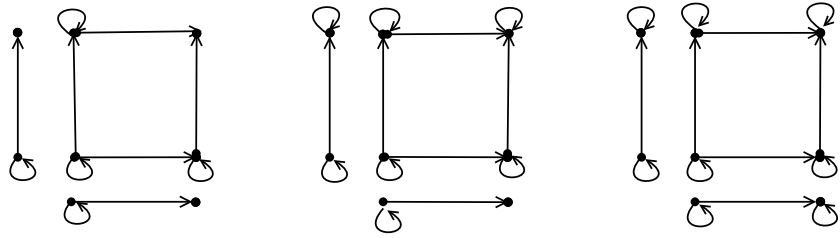}
\caption{Cartesian products of graphs with loops}
\label{fig:looped}
\end{center}
\end{figure}

For unlooped graphs the projections $p_i: G_i^v \mapsto G_i$ are isomorphism, but this does not hold for graphs with loops unless $G_i^v$ contains an unlooped vertex, see the left part of Figure \ref{fig:looped}. The other part of the figure shows that directed graphs with no unlooped vertex need not have unique prime factorizations.

\medskip

Often we color the edges or arcs of a graph and denote the color of $uv$ by $c(uv).$ For example, we usually color a product $\prod_{i=1}^kG_i$ by $k$ colors such that the edges of the $G_i$-layers are assigned color $i.$

We also need the fact that layers are {convex}, where a subgraph $H$ of a graph $G$ is called \emph{convex} if any  shortest path $P$ in $G$ between two vertices of  $H$ is already in $H.$ Another  important property we formulate as a lemma.

\begin{lemma}\label{le:projection}
 Let $G = G_1\cp G_2$ be the product of two connected graphs in $\Gamma$, and $u, v\in V(G)$. Then there exists a  unique vertex $x \in G_1^u$ of shortest distance from $v$ to  $G_1^u$, and  to any vertex $y\in G_2^v$ there is a shortest $v,y$-path that contains $x.$
\end{lemma}

\proof Let $u = (u_1,u_2)$ and $v= (v_1,v_2)$. Then  $x=(v_1,u_2)$ is the unique vertex of $G_2^v\cap G_1^u.$ By the Distance Formula \cite[Lemma 5.2]{handbook}
$$d_{G}(a, b) = d_{G_1}(a_1, b_1) + d_{G_2}(a_2, b_2)  $$
for any $a,b \in G_1\cp G_2$.
Because $x_1 = v_1$ and $x_2 = y_2$ we have  $d_G(v,x) + d_G(x,y) )=  d_{G_1}(v_1, x_1) + d_{G_2}(v_2, x_2) + d_{G_1}(x_1, y_1) + d_{G_2}(x_2,y_2)  = d_G(v,y). $ Therefore  $x$ is on a shortest $v,y$-path for any $y\in G_2^v.$

This also means that every vertex $y\in G_2^v$ that is different from $x$ must have larger distance from $v$ than $x.$ Hence $x$ is the unique vertex of shortest distance from $v$ in $G_1^u.$ \qed

It will be convenient to call the vertex $x$ of shortest distance from a vertex $v$ to a layer $G_i^u$ the \emph{projection} of $v$ into $G_i^u$ and to denote it by $p_{G_i^u}(v)$. By the projection of an arc $uv$ we mean $p_{G_i^u}(uv)=p_{G_i^u}(u)p_{G_i^u}(v)$, which may be an arc or not, and if it is an arc, the orientation is not necessarily the same as that of $uv$. If $uw \in {\mathcal{S}}(G)$ and if the restriction of the projection $p_{G_i^u}$ to $uw$ and $p_{G_i^u}(uw)$ is not an isomorphism in $G$, then we call the pair $uw$,  $p_{G_i^u}(uw)$ \emph{inconsistently directed}.

If $vu$ and $vw$ are in different layers with respect to a factorization $G_1\cp G_2$, say $vu\in  G_1^v$ and $vw\in G_2^v$, then the vertex $x = p_{G_1^w}(u) = p_{G_2^u}(w)$\footnote{Clearly $\{x\}= G_1^w\cap G_2^u$.}, together with $vuw$ induces a square $vuxw$ without diagonals. It is the only square containing $v$, $u$, and $w$ and called  \emph{product square}.
We call a  product square $vuxw$ in ${\mathcal{S}}(G_1\cp G_2)$ consistently oriented if the arc or arcs between $v$ and $u$ and between $w$ and $x$ have the same orientation as the ones between $p_{G_1^u}(v)=p_{G_1^u}(w)$ and $p_{G_1^u}(u)=p_{G_1^u}(x)$, and similarly, the arc or arcs between $v$ and $w$ and between $u$ and $x$ have the same orientation as the ones between $p_{G_2^w}(v)=p_{G_2^w}(u)$ and $p_{G_2^w}(w)=p_{G_2^w}(x)$. It means that opposite edges, say $vu$, $wx$, represent either arcs $vu$, $wx$, or arcs $uv$ and $xw$, or arcs in both directions.\bigskip

\section{Algorithms}\label{sec:alg}

In this section we present two algorithms. The first one computes the prime factorization of a connected graph  $G\in \overrightarrow{\Gamma}$ from a given prime factorization of $\mathcal{S}$$(G)$, the second computes the prime factorization of a connected graph  $G\in \overrightarrow{\Gamma_0}$, where $G$ has at least one unlooped vertex, from a given prime factorization of $\mathcal{N}$$(G)$, where ${\mathcal{N}}(G)$ denotes the graph obtained from $G$ by the removal of the loops. Both algorithms are linear in the number of arcs.

The data structures that  we use  are {incidence} and {adjacency lists}. The \emph{incidence list} of a graph in $\Gamma$ lists to every edge $e = [u,v]$ its endpoints,  whereas the \emph{adjacency list} consists of the lists of neighbors $N(v)$, $v \in V(G)$.  Every edge $e = [u,v]$ appears in $N(v)$ and $N(u)$. To both of these lists we add a pointer to the place of $e$ in the incidence list,  and in the incidence list we add  pointers to  the place of $e$ in the lists $N(v)$ and $N(u)$.

For graphs $G$ in $\overrightarrow{\Gamma}$  we form the lists for ${\mathcal{S}}(G)$ and then indicate in the incidence list for every $[u,v]$ whether it is the shadow of $uv$, $vu$, or both $uv$ and $vu$. For loops we make an entry in the adjacency list, because  $v$ is a neighbor of itself if it is looped.
Clearly the space requirement is $O(|E(G)|)$, resp. $O(|A(G)|)$ in the directed case.

We also use  a BFS-\emph{ordering} of the vertices of $G$ with respect to a root $v_0.$ It consists of the sets  $L_i$ that contain the  vertices of distance $i$ from $v_0.$ Furthermore, the vertices of $G$ are enumerated by  BFS-\emph{numbers} that satisfy the condition that BFS$(v) > $ BFS$(u)$ if the distance from $v_0$ to $v$ is larger than the one from $v_0$ to $u$. Layers through $v_0$ are called \emph{unit-layers} and vertices of unit-layers are called \emph{unit-layer vertices}.

If $vu$ is an edge or an arc, then we call $u$  a \emph{down-, cross-}, or \emph{up-neighbor} of $v$ if, respectively, $u$ is in a lower $L_i$-level than $v$, in the same, or in a higher $L_i$-level. For a given BFS-ordering we  also subdivide every list of neighbors of a vertex $v$ into lists of down-, cross- or up-neighbors. For graphs in $\Gamma$ we refer to these lists as list of \emph{down-, cross-} or \emph{up-edges}, and remark that we will not need the list of up-neighbors, resp. up-edges, in our algorithms.

If  ${\mathcal{S}}(G) = \prod_{j\in J}Z_j$  we color the edges of   ${\mathcal{S}}(G)$ as described before  with $i$ colors and subdivide every  list of down- or cross-neighbors of a vertex $v$ into sublists of different colors, according to the color of $vu$.  Furthermore, because for any subset of down- or cross-edges of a vertex $v$ with color $j$ the projection $p_{Z_j^{v_0}}(v)$  is a natural bijection into the set of down-, resp.  cross-edges of $p_{Z_j^{v_0}}(v)$ in $Z_j^{v_0}$,  we use the same order in both lists.

To describe an edge $e$ we thus need a vertex $v$ of which it is a down-, or cross edge, the fact whether it is a down-or cross-edge, its color $c(e),$ and its number in the sublist of color $c(e)$. Given the place of $v$ in the adjacency list, the color of $e$ and its sublist number, we can then find the place of $e$ in the adjacency list in  constant time. Furthermore, if we know the coordinates of $v$, then  we can find the place of $v$ in the adjacency list in time proportional to the number of coordinates in which $v$ differs from $v_0$.

By Lemma \ref{le:projection} the projection $p_{Z_i^{v_0}}(v)$ is always closer to $v_0$ than $v,$ unless $v$ already is a vertex of  $Z_i^{v_0}.$

Let $ \prod_{i=1}^kG_i$ be the prime factorization of a nontrivial, connected graph $G \in \overrightarrow{\Gamma}.$ Then  $\mathcal{S}$$(G) = \prod_{i=1}^k\mathcal{S}$$(G_i)$, where the factors $\mathcal{S}$$(G_i)$ need not be prime. Let ${\mathcal{S}}(G_i) = H_{i,1}\Box \cdots \cp H_{i, \ell(i)}$ be their prime factorizations. Then
$${\mathcal{S}} (G) = \prod_{i=1}^{k}\prod_{j=1}^{\ell(i)}H_{i,j}$$
is a representation of $\mathcal{S}$$(G)$ as a Cartesian product of prime graphs. Because the prime factorization is unique, it is the prime factorization of ${\mathcal{S}}(G)$, up to the order and isomorphisms of the factors. In other words, for any prime factorization  $\prod_{j\in J}Z_j$  of ${\mathcal{S}}(G)$, there is a partition $J =J_1\cup\cdots\cup J_k$ such that ${\mathcal{S}}(G_i) = \prod_{j\in J_i}Z_j.$

To find this partition we begin with the  partition where the $J_i$ are one-element sets. Then we combine selected sets $J_i$ until we arrive at the desired \emph{final partition}. The other partitions are called \emph{temporary partitions}.

To keep track of these operations, we create a \emph{pointer} $t_c$ with $t_c: r \rightarrow i$ if 
$r \in J_i$. We begin with the trivial partition of $J$ into one-element sets. Whenever we move from temporary partition $J_1\cup\cdots\cup J_{k}$ of $J$ to a new one, say $J'_1\cup\cdots\cup J'_{k'}$,  by combining some of the $J_i$, we update the pointers.
We also assign  the \emph{temporary color $t_c(i)$} to the edges in the $G'_i$-layer of ${\mathcal{S}}(G)$, where  ${\mathcal{S}}(G'_i)= \prod_{j\in J'_i}Z_j.$  Note that  $t_c$ produces the temporary color of any edge in constant time.

To update $t_c$
when we combine two colors, at most $|J|$ pointers have to be reset, each at constant cost. Because there are only $|J|$ colors, the total cost is $O(|J|^2)$. Recall that $|J|$ is the original number of factors. Because every vertex meets a layer of every one of the factors, we cannot have more factors than the minimum degree $\delta$ of $G$. So
$O(|J|^2) = O(\delta^2) = O(n\delta) = O(m).$

Algorithm \ref{alg:directed} specifies when colors are combined. We have prove the correctness of the algorithm and to investigate its complexity.
\medskip

\begin{algorithm}[h]
\caption{Factoring directed graphs}
\label{alg:directed}
\begin{algorithmic}
\vspace{2mm}
\REQUIRE  a connected graph $G$ in $\overrightarrow{\Gamma}\!$
\STATE \hspace{26pt}the   prime factorization $\prod_{j\in J}Z_j$ of ${\mathcal{S}}(G)$ with respect to the Cartesian product
\STATE \hspace{26pt}a \emph{BFS}-numbering of $V(G)$ with root $v_0$
   \STATE \hspace{26pt}the preceding data structure
\ENSURE the prime factors of $G_i$ of $G$
\vspace{2mm}
\STATE begin with the trivial partition of $J$
\FORALL{vertices $v$ of $G$ in their \emph{BFS} order}
\FORALL {down- and cross-edges $vu$ of $v$ in ${\mathcal{S}}(G)$}
\STATE determine its temporary color $i$
\STATE consider  the set $J_i$ of the current temporary partition of $J$
\STATE set $X = (\prod_{j\in J_i}Z_j)^{v_0}$
\STATE Project the edge $vu$ into $X$
\IF{ $vu$ and $p_X(v)p_X(u)$ are not consistently directed in $G$}
\STATE  combine the temporary colors of all down edges of $v$
\STATE scan the next vertex using the new coloring
\ENDIF
\ENDFOR
\ENDFOR
\STATE compute the products $H_i =  \prod_{j\in J_i}Z_j$  of $G,$ where $J_1\cup\cdots\cup J_k$  is the last partition of $J$  
\STATE compute the subgraphs  of $G$ induced by the $H_i$ and denote them by $G_i$
\vspace{2mm}
\end{algorithmic}
\end{algorithm}

{\bf Correctness of Algorithm \ref{alg:directed}}  \, In $L_1$ all vertices are unit-layer vertices and all edges coincide with their projection into the unit-layers of their color. Hence all pairs $vu$ and $p_X(v)p_X(u)$, where $v \in L_1$ and $vu$ a down- or cross-edge,  are consistently directed. Let this be the case for all pairs $vu$ and $p_X(v)p_X(u)$ where $v$ is in $L_0\cup \cdots \cup L_{k-1}$.

Suppose $v\in L_k$ and the algorithm detects a down- or cross-edge $vu$ of $v$ of temporary color $i$ for which $vu$ and its projection  into $X = (\prod_{j\in J_i}Z_j)^{v_0}$ are not consistently directed. Let  $Y= \prod_{j\in J\setminus J_i}Z_j$ and  $P$ be a shortest path from $v$ to $p_X(v)$. Clearly $P \in Y^v,$ and the temporary color of all edges of $P$ is different from that of $vu$.

Let $P'$ be the projection of $P$ into $Y^u.$ Then the vertices of $P$, together with those of  $P'$,  induce a subgraph $L$ of ${\mathcal{S}}(G)$ that is isomorphic to the product $P\cp [u,v]$. We call it a \emph{ladder} and  the edges of $L$ that are not in $P\cup P'$ \emph{rungs}. All rungs $ab$ project into $p_X(vu)$, and, with the exception of the pair $\{vu, p_X(vu)\}$,  all pairs $\{ab,  p_X(vu)\}$ are consistently directed, because $a$ is in some $L_j$ with $j<k$. If $v'u'$ denotes the rung closest to $vu$, then $vv'u'u$ is a product square, because $vu$ and $vv'$ have different colors. But then $vv'u'u$ is an inconsistently directed product square and we have to combine the colors of its edges.

Now, let $vq$ be any down-edge of $v$, whose temporary color $t_c(vq) \neq i$. Then there is a shortest path $Q$ from  $v$ to $p_X(v)$ that contains $q$, and hence the color of $vq$ must be combined with that of $vu$.

When the algorithm terminates, all projections are consistently oriented. We thus arrive a factorization of $G$.

We still have to show that the $G_i$ are prime, that is, that we have not merged colors unnecessarily. Notice that we have only merged colors of inconsistently oriented product squares. Originally we had the product squares of the decomposition $\prod_{j\in J}Z_j$ of ${\mathcal{S}}(G)$. Since our operations only combined colors of edges in inconsistently directed  squares,  every single combination of colors was forced.
\medskip

{\bf Complexity of Algorithm \ref{alg:directed}} \,The projection of a vertex $v$ into  $X = (\prod_{j\in J_i}Z_j)^{v_0}$ has at most $|J_i|$ coordinates in $\prod_{j\in J}Z_j$ that are different form those of $v_0$, hence $p_X(v)$ can be computed in $O(|J_i|) = O(\delta)$ time. To find the edge $p_X(vu)$ we then need its type (down- or cross-edge), its color (original color) and its sublist number, which we inherit from $vu$.  Then we can find $p_X(vu)$ in constant time and check in constant time whether the pair
$vu$, $p_X(vu)$ is consistently directed. Since the number of down-edges is bounded by $d(v)$ the time complexity for each $v$ is thus $O(\delta) + O(d(v))$, and for all vertices together it requires $O(\sum_{v \in V(G)}(\delta + d(v))) = O(m)$ time.

For the complexity of keeping track of the colorings, let us recall that we have to combine colors at most $|J|$ times and that every merging operation  of two colors costs $O(|J|)$ time, hence the overall cost is $O(|J|^2)$, of which we already know that it is $O(m)$.

\begin{theorem} Let $G$ be a connected, directed graph. Given a prime factorization of the shadow ${\mathcal{S}}(G)$ of $G$ with respect to the Cartesian product, one can compute the prime factorization of $G$ with respect to the directed product in $O(|A(G)|$ time.
\end{theorem}

\proof Given the adjacency lists of the prime factors of ${\mathcal{S}}(G)$ we can compute the data structure that we need for Algorithm \ref{alg:directed} in time and space that is linear in the number of arcs. If the factors are given, say, via their adjacency matrices, we can still find their adjacency lists in linear time and then continue in linear time and space.

If the factorization of ${\mathcal{S}}(G)$ is computed by the algorithm of Imrich and Peterin \cite{impe-2007}, then one can use the data structure provided  by that algorithm. \qed

We continue with the prime factorization of connected, directed graphs with loops.

\begin{algorithm}[h]
\caption{Factoring directed graphs with loops}
\label{alg:loops}
\begin{algorithmic}
\vspace{2mm}
\REQUIRE  a connected graph $G$ in $\overrightarrow{\Gamma_0}\!$ with an unlooped vertex
\STATE \hspace{26pt}the   prime factorization $\prod_{j\in J}Z_j$ of ${\mathcal{N}}(G)$ with respect to the Cartesian product \STATE \hspace{26pt}a \emph{BFS}-numbering of $V(G)$ with root $v_0$, which is unlooped
 \STATE \hspace{26pt}the data structure from Algorithm \ref{alg:directed}
\ENSURE the prime factors of $G_i$ of $G$.
\vspace{2mm}
\STATE begin with the trivial partition of $J$
\FORALL{vertices $v$ of $G$ in their \emph{BFS} order}
\STATE compute the projection of $v$ into the unit-layers of
the products $X_i = \prod_{j\in J_i}Z_j$, where $J_1\cup\cdots\cup J_k$ is the current temporary partition of $J$.
\IF{ $v$ is unlooped and all projections are unlooped}
\STATE  continue with the next $v$
\ELSIF{$v$ has a loop and there a projection with a loop}
\STATE  continue with the next $v$
\ELSE
\STATE{combine the temporary colors of the down edges of $v$}
\STATE continue with the next $v$
\ENDIF
\ENDFOR
\STATE compute the products $H_i =  \prod_{j\in J_i}Z_j$  of $G,$ where $J_1\cup\cdots\cup J_k$  is the last partition of $J$  
\STATE compute the subgraphs  of $G$ induced by the $H_i$ and denote them by $G_i$
\vspace{2mm}
\end{algorithmic}
\end{algorithm}

{\bf Correctness of Algorithm \ref{alg:loops}}  \, If the projection of an unlooped vertex $v$ into a unit layer (with respect to an unlooped root)  has a loop, then $v$ must be contained in that layer, and hence all shortest paths form $v$ to the projection of $v$. As in the correctness argument of Algorithm \ref{alg:directed} this means that the colors of all down-edges of $v$ have to be combined.

Suppose all projections of a vertex with a loop are unlooped. Consider all projections  that are different from $v_0$ and the respective unit-layers, that is,  the layers  $X_i^{v_0} = (\prod_{j\in J_i}Z_j)^{v_0}$, where the $J_i$ correspond to the temporary colors of the down-edges of $v$. Let $J'$ be the set of the indices of these $J_i$ and form, for any proper subset $J''$ of $J'$,  the product
$X = \prod_{j\in J''}X_j$. Then the projection of $v$ into $X^{v_0}$ is unlooped. This means, unless we combine all colors of the
down-edges of $v$ we will only have unlooped projections.

{\bf Complexity of Algorithm \ref{alg:loops}}  \,  The computation of the projection of a vertex $v$ into $(\prod_{j\in J_i}Z_j)^{v_0}$ takes $O(|J_i|)$ time. The cost of computing them all is thus $O(\sum_{i=1}^k|J_i|) = O(|J|) = O(\delta)$, which gives a total of
$O(n\delta) = O(m)$ for all vertices together.

The complexity of keeping track of the colorings is $O(m)$, as in Algorithm \ref{alg:directed}.

We have thus shown the following theorem.

\begin{theorem} Let $G$ be a connected, directed graph with loops.  Given a prime factorization of ${\mathcal{N}}(G)$ with respect to the Cartesian product, one can compute the prime factors of $G$ with respect to the Cartesian product in $O(|A(G)|$ time.
\end{theorem}\bigskip

\section{The weak Cartesian product}\label{sec:weak}

We consider infinite graphs now and begin with the definition of the Cartesian product of infinitely many factors.

\begin{definition}\label{def:prod}
Let $G_\iota$, $\iota \in I$, be a collection of graphs in $\overrightarrow{\Gamma_0}$ . Then the Cartesian product
$$G = \prod_{\iota\in I}G_\iota$$
has as its set of vertices $V(G)$ all functions
$$v: I \rightarrow \bigcup_{\iota\in I}V(G_\iota)$$
 with the property that $v(\iota)\in V(G_\iota).$ We call $v(\iota)$ the $\iota$-coordinate of $v$ and also denote it  by $v_\iota.$

The set $A(G)$ of arcs of $G$ consists of all ordered pairs $uv$ for which there exists a $\kappa$ such that $u_\kappa v_\kappa \in A(G_\kappa)$ and  $u_\iota = v_\iota$ for all $\iota \in I\setminus\{\kappa\}.$

Furthermore, $v\in V(G)$ has a loop if at least one  $v_\iota$, $\iota\in I$,  has a loop in $G_\iota.$
\end{definition}

This definition  is equivalent to the definition of the Cartesian product for two factors as given  in  Section \ref{sec:prelim}. To see this, observe that $V(G_1\cp G_2) = \{(v_1,v_2)\,|\,  v_1\in V(G_1),  v_2\in V(G_2)\}$ by the old definition.
Representing the ordered  pairs $\{(v_1,v_2)\,| \, v_1\in V(G_1),  v_2\in V(G_2)\}$  by the set of functions
$v: I \rightarrow V(G_1)\cup V(G_2)$, where $v_1\in V(G_1)$ and $v_2\in V(G_2)$,   it becomes clear  that the  definitions are equivalent.

Furthermore, let $n$ be a positive integer and  $J_1\cup\cdots\cup J_{k}$ is an arbitrary partition of the set of integers between 1 and $n$, then there clearly is a natural isomorphism between $\prod_{i=1}^nV(G_i)$, the product
$$\prod_{j=1}^{k}\prod_{i\in J_k}V(G_i),$$
and the set of functions $v: I \rightarrow \bigcup_{i=1}^nV(G_i)$, where $v_i\in V(G_i).$
This means that Cartesian multiplication of finitely many sets is associative.

For infinite $I$ and any of its partitions $\{J_\lambda\, |\, \lambda \in \Lambda\}$ we only have the isomorphism between
$\prod_{\lambda \in \Lambda} \prod_{\iota \in J_\lambda}V(G_\iota)$ and the set of functions $v: I \rightarrow \bigcup_{\iota\in I}V(G_\iota)$ where $v(\iota)\in V(G_\iota),$ but the Cartesian product of infinitely many sets is still associative.

The fact that   the set of factors need not be ordered in Definition \ref{def:prod} reflects the fact that Cartesian multiplication is commutative. 

\begin{lemma} Cartesian multiplication of directed graphs with loops is associative.
\end{lemma}

\proof Let $G_\iota$, $\iota \in I$, be a collection of directed graphs with loops, $G=\prod_{\iota\in I}G_{\iota}$, $\{J_\lambda \,|\, \lambda \in \Lambda\}$ an  arbitrary partition of $I$ and
$$H = \prod_{\lambda \in \Lambda} \prod_{\iota \in J_\lambda}G_\iota\,.$$
We  have to show that $G\cong H.$
To see this recall that there is a natural bijection between the vertices of $G$ and $H$. We will thus use the same notation for the vertices of $G$ and $H.$

Consider an arc $uv$ in $G.$ All coordinates of $u$ and $v$ are identical, except for one, say $u_\kappa \neq v_\kappa$, and $u_\kappa v_\kappa \in A(G_\kappa).$ Let $\kappa \in J_{\lambda'}$ and let $u_{J_\lambda}$ be the vertex in $\prod_{\iota \in J_\lambda}G_\iota$ with the coordinates $u_\iota$, where $\iota \in J_\lambda.$ Analogously we define
the $v_{J_\lambda}.$  Then $u_{J_\lambda} = v_{J_\lambda}$ for $\lambda \in \Lambda \setminus \{\lambda'\}$, and in this case $u_{J_{\lambda'}}v_{J_{\lambda'}}\in A(\prod_{\iota \in J_{\lambda'}}G_\iota).$
This means that $uv$ is an arc in  $H$ if it is an arc in $G.$

On the other hand, if we have an arc $uv$ in $H$, then $u_{J_\lambda} = v_{J_\lambda}$ for all $\lambda \in \Lambda$, except for one, say $\lambda'$,  for which $u_{J_{\lambda'}}v_{J_{\lambda'}}\in A(\prod_{\iota \in J_{\lambda'}}G_\iota).$ But then $u_\iota = v_\iota$ for all $\iota \in J_\lambda$ for $\lambda \neq \lambda'$, and in $J_{\lambda'}$ there is a $\kappa$ such that $u_\kappa v_\kappa \in A(G_\kappa)$ and $u_\iota v_\iota$ for $\iota \in J_{\lambda'}\setminus \{\kappa\}$, hence $uv \in A(G).$ Thus $uv$ is an arc in $G$ if it is an arc in $H.$

For the loops we observe that all $G$-coordinates of a vertex $v$ are unlooped if and only if all $H$-coordinates of $v$ are unlooped. Hence a vertex has no loop in $G$ if and only if it has no loop in $H.$
\qed

It is easily seen, and well known, that the Cartesian product of finitely many factors is connected if and only if all factors are connected, but the product of infinitely many nontrivial graphs is always disconnected. The reason is that such products contain vertices that differ in infinitely many coordinates, but every arc (which is not a loop) connects vertices that differ in exactly one coordinate, which means that the coordinates of any two vertices that are connected by a path can differ in only finitely many coordinates.

The connected components of the Cartesian product $\prod_{\iota\in I}G_\iota$ of infinitely many connected factors are called \emph{weak Cartesian products}. We will use the notation $\prod_{\iota\in I}^aG_\iota$ to indicate the component that contains the vertex $a\in V(\prod_{\iota\in I}G_\iota).$

For finite graphs, or for products of finitely many factors, the Cartesian product and the weak Cartesian product coincide.\bigskip

\section{Unique prime factorization}\label{sec:unique}

It is well known that simple, connected graphs have unique prime factorizations with respect to the weak Cartesian product; see \cite{miller-1970,im-1971}. We complete this note by showing that this is also the case for directed graphs with loops. Our proof is direct, it does not use the result for undirected graphs.

\begin{theorem} Every connected, finite or infinite graph in $\overrightarrow{\Gamma_0}$ with at least one unlooped vertex has a unique prime factor decomposition with respect to the Cartesian or weak Cartesian product.
\end{theorem}

\proof Let $G$ be a connected graph in $\overrightarrow{\Gamma_0}$ and $a$ an unlooped vertex. Consider all possible representations of $G$ as products of two factors $A_{1,\kappa}\,\cp A_{2,\kappa}$, $\kappa \in K$, and let $e$ be an arc  incident with $a.$ Observe that, for any $\kappa\in K$, exactly one of the layers ${A_{1,\kappa}^a}$, $A_{2,\kappa}^a$ contains $e.$ We denote it by $A^a_{i_e,\kappa}$ and  form
$$G_e = \bigcap_{\kappa \in K} A^a_{i_e,\kappa}\,.$$
Clearly $G_e$ is convex in $G$, because it is the intersection of convex subgraphs of $G.$ Also, if $f$ is another arc incident with $a$, and if $G_f \neq G_e$, then there exist $\kappa\in K$ such that $A^a_{i_f,\kappa}\neq A^a_{i_e,\kappa}.$ Then  $G_f\cap G_e \subseteq A^a_{i_f,\kappa}\cap A^a_{i_e,\kappa} = a.$ Hence, any two distinct
$G_f$, $ G_e$ have only the vertex $a$ in common.

We will show that  $G_e^a$  is a factor of $G$, or, to be more precise, a unit layer with basepoint $a$ of a factorization of $G$, and that $G$ is the Cartesian or weak Cartesian product of all  distinct $G_e$-s. Note, if $G_e^a$  is a factor, then it must be prime, otherwise it would not be the intersection of all layers containing $e$, and if $P$ is a prime factor of $G$, and $e$ an arc incident with $a$ that projects to an arc of $P$, then $P^a = G_e.$

Hence, if $G_\iota$, $\iota \in I$, denotes the set of all distinct $G_e$-s, and if we can show that
 $$G = \prod_{\iota\in I}^aG_\iota\,,$$
 then this is the unique representation of $G$ as a Cartesian or weak Cartesian product.

We now define the coordinates of a vertex $v$ of $G$ as the vertices of shortest distance in the $G_e$-s from $v.$ So the $G_e$-coordinate of $v\in G$ is the unique vertex of $G_e$ that is closest to $v$. Clearly all coordinates of $a$ are $a.$

Consider a vertex  $v \neq a.$ If $v$ is either in $A_{1,\kappa}^ a$ or in $A_{2,\kappa}^a$ for all $\kappa \in K$, then $v$ is in $G_e$, where $e$ is an arc incident with $a$ on an arbitrary shortest path from $a$ to $v.$ So it is identical to its  coordinate in $G_e.$ All other coordinates clearly are $a.$
Otherwise there must be a $\kappa_1\in K$ such that $v \not\in A_{i,\kappa_1}^a$, for $i \in \{1,2\}.$
Let  $v_1$ be the projection of $v$ into ${A_{i_1,\kappa_1}^a}$, where $i_1$ is chosen from $\{1,2\}.$ Notice that $d_G(v_1,a) < d_G(v,a).$ If there is a $\kappa_2$ such that $v_1 \not\in A_{i,\kappa_2}^a$, for  $i \in \{1,2\}$, then we form $v_2$ as the projection of $v_1$ into ${A_{i_2,\kappa_2}^a}$, $i_2\in \{1,2\}.$ As before
$d_G(v_2,a) < d_G(v_1,a)$, hence this process cannot continue indefinitely. It ends when, for some $k$,  $v_k = a$, or when $v_k$ is either in $A_{1,\kappa}^ a$ or in $A_{2,\kappa}^a$ for all $\kappa \in K$, and hence, by the previous argument, in a $G_e$, where $e$ is an arc incident with $a$ that is on a shortest path to $v_k.$
Notice that $v \notin A_{i_1,\kappa_1}^a$ and that
$$A_{i_1,\kappa_1}^a\supset A_{i_2,\kappa_2}^a \supset \cdots \supset   A_{i_k,\kappa_k}^a \supseteq G_e.$$
Hence, by Lemma \ref{le:projection}, for any vertex $x\in G_e$ there is a shortest $v,x$-path that passes through $v_1.$ Similarly, there is a shortest
$v_1,x$-path through $v_2$, and finally through $v_k.$ Hence, if $x$ is on a shortest path from $G_e$ to $v$, then it must be equal to $v_k.$

Thus, $v_k$  is the unique vertex of shortest distance from $G_e$ to $v$, and thus its $G_e$-coordinate.  Note that,  by the construction of $v_k$,  every $v$ can have at most finitely many coordinates that are different from $a.$

We still have to show that different vertices have different coordinates. Let $u \neq v$, where $u\neq a$ and $e$ is an arc incident with $a$ on a shortest $a,u$-path. Suppose first that, for any $\kappa \in K$, the vertices $u$ and $v$ are either both in $A_{1,\kappa}^a$, or both in $A_{2,\kappa}^a.$ Then $u$ and $v$ are both in $G_e$, and equal to their $G_e$ coordinates, which are different. If this is not the case, then there is a $\kappa\in K$, such that the projections of  $u$ and $v$ into at least one of $A_{1,\kappa}^a$ or  $ A_{2,\kappa}^a$ are distinct and different from $\{u,v\}$, we call them $u_1,v_1.$  We can continue this until both $u_k$ and $v_k$ are in the same $G_e.$ But then the $G_e$-coordinate of $v$ is $a$ and different from the $G_e$-coordinate of $u.$

Having uniquely coordinatized the vertices of $G$, let us consider the projections of loops and arcs into the factors.

First the arcs. Let $uv$ be an arc in $G.$ If $uv \in G_e$, then $uv$ is equal to its projection into $G_e$, and all the other projections, that is, coordinates, are $a.$ Otherwise we proceed as above, when we showed that different vertices have different coordinates. We just have to observe that $uv$ is an arc if and only if $u_1v_1$ is an arc, and, by induction, that $uv$ is an arc if and only if the projections into $G_e$ are an arc. Furthermore, there is always only one projection we can choose from such that the $u_i,v_i$ remain different. So, if we alter any of the  sequences of projections that sends  $u$ into $G_e$, then $u$, $v$ end up in the same vertex in some other $G_f$, which means that they have the same $G_f$-coordinates.

Now the loops. Observe, if $v$ is unlooped, then this is also the case of $v_1, v_2, \ldots, v_k$, but if $v$ has a loop then one or both of the projections into ${A_{1,\kappa_1}^a}(v)$ or ${A_{2,\kappa_1}^a}(v)$ have a loop. We choose for $v_1$ a projection with  a loop and continue like this. Then this process ends in a  coordinate of $v$ that has  a loop. In other words,  $v$ is unlooped if and only if all projections are unlooped.
\qed


\end{document}